\documentclass[12pt,reqno]{amsart}
\usepackage{amsmath}
\usepackage{amsxtra}
\usepackage{amscd}
\usepackage{amsthm}
\usepackage{amsfonts}
\usepackage{amssymb}
\usepackage{eucal}

\newtheorem{theorem}{Theorem}[section]
\newtheorem{conj}[theorem]{Conjecture}

\newtheorem{lem}[theorem]{Lemma}

\theoremstyle{definition}

\theoremstyle{remark}
\newtheorem{rem}[theorem]{Remark}
\theoremstyle{remark}

\numberwithin{equation}{section}

\newcommand{\nc}{\newcommand}
\nc{\on}{\operatorname}
\nc{\bs}{\boldsymbol}
\nc{\ri}{\rangle}
\nc{\lef}{\langle}
\nc{\Res}{\on{Res}}
\nc{\mc}{\mathcal}
\nc{\mbf}{\mathbf}
\nc{\bb}{{\mathfrak b}}

\nc{\Z}{{\mathbb Z}}
\nc{\C}{{\mathbb C}}
\nc{\R}{{\mathbb R}}
\nc{\Q}{{\mathbb Q}}

\nc{\g}{{\mathfrak g}}
\nc{\h}{{\mathfrak h}}
\nc{\n}{{\mathfrak n}}

\nc{\al}{\alpha}
\nc{\la}{\lambda}
\nc{\La}{\Lambda}
\nc{\ep}{\epsilon}
\nc{\om}{\omega}
\nc{\Om}{\Omega}
\nc{\ga}{\gamma}
\nc{\Ga}{\Gamma}

\nc{\bean}{\begin{eqnarray}}
\nc{\eean}{\end{eqnarray}}
\nc{\be}{\begin{displaymath}}
\nc{\ee}{\end{displaymath}}
\nc{\bea}{\begin{eqnarray*}}   
\nc{\eea}{\end{eqnarray*}}

\nc{\bi}{\bibitem}
\nc{\Ref}[1]{{$($\ref{#1}$)$}}

\begin{document}

\title[Alternating sums]
{Factorization of alternating sums of Virasoro characters}
\author[E. Mukhin]
{E. Mukhin}
\address{Department of Mathematical Sciences, Indiana University -
Purdue University Indianapolis, 402 North Blackford St, Indianapolis,
IN 46202-3216, USA, \newline mukhin@math.iupui.edu}

\maketitle

\begin{abstract} 
  G. Andrews proved that if $n$ is a prime number then
  the coefficients $a_k$ and
  $a_{k+n}$ of the product $(q,q)_\infty/(q^n,q^n)_\infty=\sum_k a_kq^k$
  have the same sign, see \cite{A}. 
  We generalize this result in several directions.
  Our results are based on the observation that many products can be
  written as alternating sums of characters of Virasoro modules.
\end{abstract}

\section{Introduction} 
In the past several decades with the appearance and rapid development
of Conformal Field Theory, the Virasoro modules enjoyed ample
attention from mathematicians and physicists. 

In this paper we study the following question.

\medskip 
{\bf Question.} Which finite (alternating) sum of characters of 
Virasoro modules occurring in a minimal series can be written 
in the form 
$$
\frac{\prod_{i\in I_-}(1-q^i)\ \prod_{i\in I_+}(1+q^i) }{\prod_{j\in J_-}(1- q^j)\ \prod_{j\in J_+}(1+ q^j)}\ ,
$$
where $I_\pm,J_\pm$ are some sets of natural numbers?

\medskip

This  question is motivated by a number of applications to combinatorics 
and mathematical physics. 
 \begin{itemize} 
\medskip  
 \item The character of any single Virasoro module occurring in $(2,2r+1)$
  minimal series is factorizable. The same character can be written in a
  so-called fermionic form and we obtain the celebrated
  Rogers-Ramanujan-Gordon-Andrews identities, see \cite{A1}.\linebreak

\item For any $1<s<p'$, the 
sums and differences of characters of $(1,s)$ and $(p-1,s)$ 
modules in $(p,p')$ minimal series, where $p\in\{3,4\}$, are factorizable. 
For the case of the sum we also have a fermionic formula and 
therefore an identity of Rogers-Ramanujan-Gordon-Andrews type, see \cite{FFW}.
\medskip 
\item A sum of three Virasoro characters from $(2,9)$ minimal series equals the product $(q,q)_\infty/(q^3,q^3)_\infty$. 
It immediately implies that if $(q,q)_\infty/(q^3,q^3)_\infty =\sum_k a_k q^k$
and if $a_k$ and  $a_{k+3}$ are non zero
then they have the same sign. The famous Borwein conjecture 
is a finitization of this fact: it asserts that if coefficients $a_k^{(N)}$ and
 $a_{k+3}^{(N)}$ of the product $\prod_{j=0}^N(1-q^{3j+1})(1-q^{3j+2})=
  \sum_k a_k^{(N)} q^k$ are non-zero then they have the same sign, see
  \cite{A}.  
\medskip 
\item The number $1$ can be written in several ways as an alternating
sum of characters of Virasoro modules from $(2,2r+1)$ minimal series. 
Each such equality gives a family of non-trivial partition identities, see 
\cite{MMO}. 
\medskip   
\item Factorized form of graded characters of Virasoro modules is
  crucial for studying form factors of integrable deformations of
  Conformal Field Theory, see \cite{Ch}.
\medskip 
\item Writing a sum of graded characters of Virasoro modules in a
  product form leads to identities involving sums of products of
  graded characters of Virasoro modules.
\end{itemize}
\medskip

We present two large families of products which are equal to finite
alternating sums of Virasoro characters. 
In particular, these families contain all the known 
examples of such phenomena.

We prove our formulae by an application of the triple Jacobi
identity and of the quintuple identity to the Rocha-Caridi
formula for the characters of Virasoro modules.

In the case of the triple Jacobi identity we prove the following formula
\bean\label{for}
\frac{(q^{\frac{B(a'-c)}{2}},q^{\frac{B(a'+c)}{2}},q^{Ba'};q^{Ba'})_\infty}
{(q^n;q^n)_\infty}=q^{\frac{(p-p')^2-(cB)^2}{8Ba'}}\hspace{-20pt}
\sum_{ 0<r<p/b, \ r\equiv 1\ (\on{mod}\  2), \atop 0<s<p'/b', 
\ ps\equiv bc\ (\on{mod}\ a')} \hspace{-20pt}(-1)^{t_{r,s}} 
\chi_{rb,sb'}^{(p,p')}(q^n)\ ,
\eean
where
\be
t_{r,s}=\frac{p'r/b'-ps/b+c}{2}\ .
\ee
Here $p,p'$ are two relatively prime positive integers, and
$a',b,b',c$ are positive integers such that $2b$ divides $p$, $a'b'$
divides $p'$, $a'>c$ and $c$ is odd.  The numbers $B$ and $n$ are
given by $B=bb'$, $n=\frac{pp'}{2a'bb'}$. The right hand side of our
formula contains an alternating sum of $n$ different Virasoro
characters from $(p,p')$ minimal series. 

A few cases of such formula are known. 
The cases of  $n=1,2$ can be found in \cite{BF}, \cite{FFW}. 
The case of $b=b'=B=1$, $a'=3n$, $n=c$ can be found in \cite{MMO} 
(in this case the left hand side clearly equals to $1$).

If $n$ is even (that is if $p$ is divisible by $4$), the signs 
in the formula can be written in a different way. 
Namely, for the case of even $n$ we also have:
\be
\frac{(q^{\frac{B(a'-c)}{2}},-q^{\frac{B(a'+c)}{2}},-q^{Ba'};-q^{Ba'})_\infty}
{(q^n;q^n)_\infty}=q^{\frac{(p-p')^2-(cB)^2}{8Ba'}}\hspace{-20pt}
\sum_{ 0<r<p/b, \ r\equiv 1\ (\on{mod}\  2), \atop 0<s<p'/b', 
\ ps\equiv bc\ (\on{mod}\ a')} \hspace{-20pt}(-1)^{\frac{t_{r,s}(t_{r,s}+1)}{2}} 
\chi_{rb,sb'}^{(p,p')}(q^n)\ .
\ee
if $a'-c$ is divisible by $4$, and 
\be
\frac{(-q^{\frac{B(a'-c)}{2}},q^{\frac{B(a'+c)}{2}},-q^{Ba'};-q^{Ba'})_\infty}
{(q^n;q^n)_\infty}=q^{\frac{(p-p')^2-(cB)^2}{8Ba'}}\hspace{-20pt}
\sum_{ 0<r<p/b, \ r\equiv 1\ (\on{mod}\  2), \atop 0<s<p'/b', 
\ ps\equiv bc\ (\on{mod}\ a')} \hspace{-20pt}(-1)^{\frac{t_{r,s}(t_{r,s}-1)}{2}} 
\chi_{rb,sb'}^{(p,p')}(q^n)\ .
\ee
if $a'-c$ is not divisible by $4$. Some cases of these formulae with $n=2$ are 
contained in \cite{BF}, \cite{FFW}.

The above formulae can be used to obtain identities of
Rogers-Ramanujan-Gordon-Andrews type by equating the product side with
any known expression for Virasoro characters in the right hand side.
In particular, we expect that the known fermionic expressions for the
Virasoro characters appearing in our formula (see e.g. \cite{BM},
\cite{BMS}, \cite{Wl}) can be summed up to a fermionic form. 
We do not discuss fermionic formulae in this paper.

Another set of identities is obtained by multiplying or dividing the
product forms for different cases of the above formula.  Such
identities involve sums of products of Virasoro characters, see Remark
\ref{rem}.

\medskip

We use our formulae to study signs of coefficients of products.
Write the left hand side of \Ref{for} as a formal power series
$\sum_{j=0}^\infty \phi_jq^j$.  We conjecture that $\phi_j$ and
$\phi_{j+n}$ always have the same sign and prove it in several cases,
see Theorem \ref{pos thm}.  The case $a'=3$, $B=c=1$ and prime $n$ was
proved in \cite{A} in relation to the Borwein conjecture.

In some cases (e.g. when $B=1$ and all odd prime divisors of $n$
divide $a'$) for each $j$ there is only one term on the right hand side
of \Ref{for} which has a non-trivial coefficient of $q^j$, and it
follows that $\phi_j$ and $\phi_{j+n}$ do have the same sign.  In more
complicated cases, one can hope to make use of some fermionic
expressions for Virasoro characters to perform the subtraction.  We
use this idea to prove our conjecture for the case of prime $n$ and odd $B$.

\medskip

The formulae and results in the case of the quintuple identity are
similar, see formulae \Ref{main5 for},  \Ref{main5 for a},  
\Ref{main5 for b},  \Ref{main5 for c} and Theorem \ref{pos thm 5}.

Our paper is structured as follows.  We recall basic facts about
Virasoro modules in Section \ref{mm sec}.  Sections \ref{3p sec} and
\ref{5p sec} contain statements of the main results in the cases of
the triple and quintuple products respectively. The proofs are
collected in Section \ref{proof sec}.

\section{Minimal models}\label{mm sec}
Let ${\mc Vir}$ be the Virasoro algebra with the standard $\C$-basis
$\{L_n\}_{n\in\Z}$ and $C$, satisfying
\be
[L_m,L_n]=(m-n)L_{m+n}+\frac{C}{12}m(m^2-1)\delta_{m+n,0}, \quad
[C,L_n]=0.  
\ee 
Let $(p,p')$ be a pair of relatively prime integers greater than 1.
There exists a family of irreducible ${\mc Vir}$-modules
$M_{r,s}^{(p,p')}$ where $1\le r\le p-1$, $1\le s\le p'-1$ on which
$C$ acts as the scalar
\be
C_{p,p'}=1-\frac{6(p-p')^2}{4pp'}\ .  
\ee 
The module $M^{(p,p')}_{r,s}$ is $\Q$-graded with respect to the
degree operator $L_0$ and the corresponding formal character
$\chi_{r,s}^{(p,p')}(q):=\on {Tr}(q^{L_0})$ is given by the following
bosonic formula, see \cite{RC}:
\be
\chi_{r,s}^{(p,p')}(q):=\frac{q^{\Delta_{r,s}^{(p,p')}}}{(q)_\infty}
\left(\sum_{j\in\Z}q^{pp'j^2+(p'r-ps)j}-
  \sum_{j\in\Z}q^{pp'j^2+(p'r+ps)j+rs}\right).  
\ee
Here $(q)_\infty=\prod_{j=1}^\infty(1-q^j)$ and the conformal
dimension $\Delta_{r,s}^{(p,p')}$ is given by
\be
\Delta_{r,s}=\frac{(p'r-sp)^2-(p'-p)^2}{4pp'}\ .  
\ee 
It is convenient to write $\chi_{r,s}^{(p,p')}(q)$ in the following form:
\bean\label{RC}
\chi_{r,s}^{(p,p')}(q)=\frac{q^{-\frac{(p'-p)^2}{4pp'}}}{(q)_\infty}
\left(\sum_{j\in\Z}
q^{\frac{(2pp'j+p'r-ps)^2}{4pp'}} -
\sum_{j\in\Z}q^{\frac{(2pp'j+p'r+ps)^2}{4pp'}}\right).
\eean

From  formula \Ref{RC} one immediately observes that
\bean\label{dual}
\chi_{r,s}^{(p,p')}(q)=\chi_{p-r,p'-s}^{(p,p')}(q).
\eean

The normalized character $\bar
\chi_{r,s}^{(p,p')}$ given by
\be
\bar\chi_{r,s}^{(p,p')}(q):=q^{-\Delta_{r,s}^{(p,p')}}
\chi_{r,s}^{(p,p')}(q)=1+o(1)
\ee
is a formal power series in $q$ with non-negative coefficients. (In
fact the only zero coefficient is the coefficient of $q$ in $\bar
\chi_{1,1}^{(p,p')}=\bar\chi_{p-1,p'-1}^{(p,p')}$.)

\section{Main results}

For integers $a,b$, we write $a\perp b$ if $a,b$ are relatively prime.

Let $p,p'$ be relatively prime integers greater than $1$.  Let
$a,a',b,b'$ be natural numbers and $c$ a non-negative integer
such that $a'>c$, $ab$ divides $p$,
$a'b'$ divides $p'$.
 
We call numbers $b,b'$ {\it the scaling factors}, numbers $a,a'$ {\it
  the moduli} and number $c$ {\it the common residue}.  We obviously
have $ab \perp a'b'$.

We define
\be
B:=bb', \qquad n:=\frac{pp'}{aa'bb'}.
\ee

We use the notation $(u_1,\dots,u_k;v)_\infty:=\prod_{i=0}^\infty
\prod_{j=1}^k(1-u_jv^i)$.

\subsection{Triple products} \label{3p sec}
In this section we assume that $c$ is odd and set 
\be
a=2.
\ee
Then $p',c,a',b'$ are all odd, $p/b$ is even.

We call a pair of integers $(r,s)$ {\it 2-contributing of the first
  type} if 
\be 0< r< p/b, \qquad 0< s < p'/b', \qquad
\frac{p'r/b'-ps/b+c}{2aa'}\in\Z.  
\ee 
We call a pair of integers
$(r,s)$ {\it 2-contributing of the second type} if 
\be 0< r< p/b,\qquad
0< s < p'/b', \qquad \frac{p'r/b'+ps/b-c}{2 aa'}\in\Z.  
\ee 
We denote
the set of all 2-contributing pairs of type $j$ by $\mc A^{(2)}_j$,
$j=1,2$. We call pair of integers $(r,s)$ {\it 2-contributing} if
$(r,s)$ is either 2-contributing of the first type or 2-contributing
of the second type. We denote
the set of all 2-contributing pairs by $\mc A^{(2)}$.

\begin{lem}
We have $\mc A_1^{(2)}\bigcap \mc A_2^{(2)}=\emptyset$.  

If $\ 0<r<r+a<p/b\ $, then $(r,s)\in\mc A_j^{(2)}$ if and only if
$(r+a,s)\in\mc A_{3-j}^{(2)}$.

If $\ 0<s<s+a'<p'/b'\ $ and $p$ is even 
then $(r,s)\in\mc A_j^{(2)}$ if and only if
$(r,s+a')\in\mc A_j^{(2)}$. 

If $\ 0<s<s+a'<p'/b'\ $ and $p$ is odd
then $(r,s)\in\mc A_j^{(2)}$ if and only if
$(r,s+a')\in\mc A_{3-j}^{(2)}$.

\hspace{14cm}$\Box$
\end{lem}

\begin{theorem}\label{main} 
We have the following identity of formal power series in $q$:
\bean\label{main for}
\frac{(q^{\frac{B(a'-c)}{2}},q^{\frac{B(a'+c)}{2}},q^{Ba'};q^{Ba'})_\infty}
{(q^n;q^n)_\infty}=
q^{\frac{(p-p')^2-(cB)^2}{4Baa'}}\Big(\hspace{-10pt}\sum_{(r,s)\in \mc A_1^{(2)}}
\chi_{rb,sb'}^{(p,p')}(q^n)-\hspace{-10pt}\sum_{(r,s)\in \mc A_2^{(2)}}
\chi_{rb,sb'}^{(p,p')}(q^n)\Big) .
\eean
\end{theorem}
Theorem \ref{main} is proved in Section \ref{Proof of main}.  

The cases $n=1,2$ of Theorem \ref{main} can be found in \cite{BF}, see
also \cite{FFW}. The case   $b=b'=B=1$, $a'=3n$, $n=c$ of Theorem
\ref{main} can be found in \cite{MMO}.

We note that there are $n$ summands on the left hand side of \Ref{main
  for}, moreover, $\chi_{rb,sb'}^{(p,p')}$ is present only if
$\chi_{p-rb,p'-sb'}^{(p,p')}$ is not present, see Lemma \ref{adm}.  We
also note that formula \Ref{main for} remains the same if $c$ is
changed to $-c$, the right hand side for the obvious reason and the
left hand side because of relation \Ref{dual}.

If $n$ is even then there is a formula which differs from \Ref{main for} 
only by the choice of signs.

For a 2-contributing pair $(r,s)$ we define the integer $t_{r,s}$ 
by the formula
\be
t_{r,s}=(p'r/b'-ps/b+c)/2.
\ee
The integer  $t_{r,s}$ is even if $(r,s)\in \mc A_1^{(2)}$ and odd if 
$(r,s)\in \mc A_2^{(2)}$.

\begin{theorem}\label{main a} 
Let $n$ be even and let $a'\equiv c\ (\on{mod}\ 4)$. 
Then we have the following identity of formal power series in $q$:
\bean\label{main for a}
\frac{(q^{\frac{B(a'-c)}{2}},-q^{\frac{B(a'+c)}{2}},-q^{Ba'};-q^{Ba'})_\infty}
{(q^n;q^n)_\infty}=
q^{\frac{(p-p')^2-(cB)^2}{4Baa'}}
\Big(\hspace{-10pt}\sum_{(r,s)\in \mc A^{(2)}}
(-1)^{\frac{t_{r,s}(t_{r,s}+1)}{2}}\chi_{rb,sb'}^{(p,p')}(q^n)\Big)
.
\eean

Let $n$ be even and let $a'\not\equiv c\ (\on{mod}\ 4)$.
Then we have the following identity of formal power series in $q$:
\bean\label{main for b}
\frac{(-q^{\frac{B(a'-c)}{2}},q^{\frac{B(a'+c)}{2}},-q^{Ba'};-q^{Ba'})_\infty}
{(q^n;q^n)_\infty}=
q^{\frac{(p-p')^2-(cB)^2}{4Baa'}}
\Big(\hspace{-10pt}\sum_{(r,s)\in \mc A^{(2)}}
(-1)^{\frac{t_{r,s}(t_{r,s}-1)}{2}}\chi_{rb,sb'}^{(p,p')}(q^n)\Big)
.
\eean
\end{theorem}
Theorem \ref{main a} is proved in Section \ref{Proof of main a}.  

Some cases with $n=2$ of Theorem \ref{main a} can be found in \cite{BF}, see
also \cite{FFW}.
\medskip

We apply Theorem \ref{main}
to study the signs of the coefficients 
of products. 

Fix natural numbers $a',B, c,n$ such that $a'>c$, $a'c\perp 2$.
Define formal power series $\phi(q)$  by the formula:
\be
\phi(q)=\phi_{a',B, c,n}(q):=
\frac{(q^{\frac{B(a'-c)}{2}},q^{\frac{B(a'+c)}{2}},q^{Ba'};q^{Ba'})_\infty}
{(q^n;q^n)_\infty}\ .
\ee

We note that $\phi_{ka',B, kc,n}(q)=\phi_{a',kB, c,n}(q)$ and
$\phi_{a',kB, c,kn}(q)=\phi_{a',B,c,n}(q^k)$. 
Therefore without loss of generality we assume $a'\perp c$ and $B\perp n$.

We write
\be 
\phi(q)=\sum_{j=0}^\infty \phi_{j}q^j.
\ee

\begin{conj}\label{conj}
We have $\phi_{j}\phi_{j+n}^{(i)}\geq 0$ for all $j\in\Z_{\geq 0}$.
\end{conj}

If $n=1$ then all factors in the numerator of $\phi(q)$
cancel with factors in the denominator and therefore all coefficients
$\phi_j$ are clearly positive. 

By Theorem \ref{main}, we can always write $\phi(q)$ as a
sum of Virasoro characters (usually in several ways).  This fact can
be used to prove several cases of Conjecture \ref{conj}.

\begin{theorem}\label{pos thm} 
Conjecture \ref{conj} holds in each of the following cases:
\begin{enumerate}
\item all odd prime divisors of $n$ divide $a'$;
\item $n$ is a prime number, $B$ is odd.
\end{enumerate}
\end{theorem}

Theorem \ref{pos thm} is proved in Section \ref{Proof of pos}.
Theorem \ref{pos thm} in the case of $a'=3,c=1,B=1$ and prime $n$ is
proved in \cite{A}.

\begin{rem}
  It follows immediately from Theorem \ref{main} that coefficient
  $\phi_j$ is zero unless there exists a 2-contributing pair $(r,s)$
  such that $((p'rb-psb')^2-(cB)^2)/(4Baa')-j$ is divisible by $n$.
  Equivalently, coefficient $\phi_j$ is zero unless there exists an
  integer $m$ such that $mB(a'm+c)/2-j$ is divisible by $n$, see Section
  \ref{Proof of main}.
\end{rem}

\begin{rem}
The results similar to Theorem \ref{pos thm} also hold for 
the products appearing in the left hand sides of
formulae \Ref{main for a} and \Ref{main for b}.
\end{rem}

\subsection{Quintuple products}\label{5p sec}
We set
\be
a=3.
\ee
Then $p',a',b'$ are not divisible by 3, 
$p$ is divisible by 3.

We call a pair of integers $(r,s)$ {\it 3-contributing of the first
  type} if 
\be 0< r< p/b, \qquad 0< s < p'/b', \qquad
\frac{p'r/b'-ps/b-a'+3c}{2aa'}\in\Z.  
\ee 
We call a pair of integers
$(r,s)$ {\it 3-contributing of the second type} 
\be 0< r< p/b,\qquad
0< s < p'/b', \qquad \frac{p'r/b'+ps/b+a'-3c}{2aa'}\in\Z.  
\ee 
We denote the set of all 3-contributing pairs of type $j$ by $\mc
A_j^{(3)}$, $j=1,2$. We call pair of integers $(r,s)$ {\it 3-contributing} if
$(r,s)$ is either 3-contributing of the first type or 3-contributing
of the second type. We denote
the set of all 3-contributing pairs by $\mc A^{(3)}$.

\begin{lem}
We have $\mc A_1^{(3)}\bigcap \mc A_2^{(3)}=\emptyset$.  

If $\ 0<r<r+2a<p/b\ $, then $(r,s)\in\mc A_j^{(3)}$ if and only if
$(r+2a,s)\in\mc A_j^{(3)}$. 

If $\ 0<s<s+2a'<p'/b'\ $, then $(r,s)\in\mc
A_j^{(3)}$ if and only if $(r,s+2a')\in\mc A_j^{(3)}$. \qquad $\Box$

\end{lem}

\begin{theorem}\label{main5} 
We have the identity of formal power series in $q$:
\bean\label{main5 for}
\frac{(q^{Bc},q^{B(2a'-c)},q^{2Ba'};q^{2Ba'})_\infty 
(q^{2B(a'+c)},q^{2B(a'-c)};q^{4Ba'})_\infty}{(q^n;q^n)_\infty}
 =\hspace{120pt}\notag\\
q^{\frac{(p-p')^2-(a'-3c)^2B^2}{4Baa'}}
\left(\sum_{(r,s)\in \mc A_1^{(3)}}
\chi_{rb,sb'}^{(p,p')}(q^n)- \sum_{(r,s)\in \mc A_2^{(3)}}
\chi_{rb,sb'}^{(p,p')}(q^n)\right) .
\eean
\end{theorem}
Theorem \ref{main5} is proved in Section \ref{Proof of main5}.  The
cases $n=1,2$ of Theorem \ref{main5} can be found in \cite{BF}, see
also \cite{FFW}.

We note that there are $n$ summands on the left hand side of formula
\Ref{main5 for}, moreover, $\chi_{rb,sb'}^{(p,p')}$ is present only if
$\chi_{p-rb,p'-sb'}^{(p,p')}$ is not present, see Lemma \ref{adm5}.

If $n$ is even then we have formulae which differ from \Ref{main5 for} 
only by the choice of signs.

For a 3-contributing pair $(r,s)$ define the integer $f_{r,s}$ as follows.
If $(r,s)$ is a 3-contributing pair of the first kind we set
\be
f_{r,s}=\frac{p'r/b'-ps-a'+3c}{2aa'}.
\ee
If $(r,s)$ is a 3-contributing pair of the second kind we set
\be
f_{r,s}=\frac{p'r/b'+ps-a'+3c}{2aa'}.
\ee

\begin{theorem}\label{main5 a} 
Let $p'/b'$ be even or let $p/b$ be even and $c$ odd.
We have the identity of formal power series in $q$:
\bean\label{main5 for a}
\frac{(-q^{Bc},-q^{B(2a'-c)},q^{2Ba'};q^{2Ba'})_\infty 
(q^{2B(a'+c)},q^{2B(a'-c)};q^{4Ba'})_\infty}{(q^n;q^n)_\infty}
 =\hspace{100pt}\notag\\
q^{\frac{(p-p')^2-(a'-3c)^2B^2}{4Baa'}}
\left(\sum_{(r,s)\in \mc A_1^{(3)}}(-1)^{f_{r,s}}
\chi_{rb,sb'}^{(p,p')}(q^n)- \sum_{(r,s)\in \mc A_2^{(3)}}
(-1)^{f_{r,s}}
\chi_{rb,sb'}^{(p,p')}(q^n)\right) .
\eean

Let $p'/b'$ be divisible by $4$ or let $p/b$ and $c$ be divisible by $4$.
We have the identity of formal power series in $q$:
\bean\label{main5 for b}
\frac{(q^{Bc},-q^{B(2a'-c)},-q^{2Ba'};-q^{2Ba'})_\infty 
(-q^{2B(a'+c)},-q^{2B(a'-c)};q^{4Ba'})_\infty}{(q^n;q^n)_\infty}
 =\hspace{70pt}\notag\\ 
q^{\frac{(p-p')^2-(a'-3c)^2B^2}{4Baa'}}
\left(\sum_{(r,s)\in \mc A_1^{(3)}}\hspace{-10pt}
(-1)^{\frac{f_{r,s}(f_{r,s}-1)}{2}}
\chi_{rb,sb'}^{(p,p')}(q^n)- \hspace{-10pt}
 \sum_{(r,s)\in \mc A_2^{(3)}} \hspace{-10pt}
(-1)^{\frac{f_{r,s}(f_{r,s}-1)}{2}}
\chi_{rb,sb'}^{(p,p')}(q^n)\right) .
\eean

Let $p'/b'$ be divisible by $4$ or let $p/b$ and $c+2$ be divisible by $4$.
We have the identity of formal power series in $q$:
\bean\label{main5 for c}
\frac{(-q^{Bc},q^{B(2a'-c)},-q^{2Ba'};-q^{2Ba'})_\infty 
(-q^{2B(a'+c)},-q^{2B(a'-c)};q^{4Ba'})_\infty}{(q^n;q^n)_\infty}
 =\hspace{70pt}\notag\\
q^{\frac{(p-p')^2-(a'-3c)^2B^2}{4Baa'}}
\left(\sum_{(r,s)\in \mc A_1^{(3)}}\hspace{-10pt}
(-1)^{\frac{f_{r,s}(f_{r,s}+1)}{2}}
\chi_{rb,sb'}^{(p,p')}(q^n)- \hspace{-10pt}
\sum_{(r,s)\in \mc A_2^{(3)}}
\hspace{-10pt}(-1)^{\frac{f_{r,s}(f_{r,s}+1)}{2}}
\chi_{rb,sb'}^{(p,p')}(q^n)\right) .
\eean
\end{theorem}
Theorem \ref{main5 a} is proved in Section \ref{Proof of main5 a}. 

We apply Theorem \ref{main5} to study the signs of the coefficients 
of products. 

Fix natural numbers $a',B, c,n$ such that $a'>c$, $a'\perp 3$.  Define
the formal power series $\psi(q)$ by the formula:
\be
\psi(q)=\psi_{a',B, c,n}(q):=
\frac{(q^{Bc},q^{B(2a'-c)},q^{2Ba'};q^{2Ba'})_\infty 
(q^{2B(a'+c)},q^{2B(a'-c)};q^{4Ba'})_\infty}{(q^n;q^n)_\infty}\ .
\ee
We note that $\psi_{ka',B, kc,n}(q)=\psi_{a',kB, c,n}(q)$ and
$\psi_{a',kB, c,kn}(q)=\psi_{a',B,c,n}(q^k)$. Therefore without loss
of generality we assume $a'\perp c$ and $B\perp n$.

We write
\be 
\psi(q)=\sum_{j=0}^\infty \psi_{j}q^j.
\ee

\begin{conj}\label{conj5}
We have $\psi_{j}\psi_{j+n}\geq 0$ for all $j\in\Z_{\geq 0}$.
\end{conj}
If $n=1$ then all factors in the numerator of $\psi(q)$ cancel with
factors in the denominator and therefore all coefficients $\psi_j$ are
clearly positive. Thus Conjecture \ref{conj} is obviously true when
$n=1$.

Note that by Theorem \ref{main}, we can always write $\psi(q)$ as a
sum of Virasoro characters (usually in several ways).  This fact can
be used to prove some cases of Conjecture \ref{conj5}.

\begin{theorem}\label{pos thm 5} 
  Conjecture \ref{conj5} is true if all prime divisors of $n$ different
  from 3 divide $a'$.
\end{theorem}

Theorem \ref{pos thm 5} is proved in Section \ref{Proof of pos 5}.

\begin{rem}
The results similar to Theorem \ref{pos thm 5} also hold for 
the products appearing in the left hand sides of
formulae \Ref{main5 for a},\Ref{main5 for b} and \Ref{main5 for c}.
\end{rem}

\begin{rem}
  It follows immediately from Theorem \ref{main5} that coefficient
  $\psi_j$ is zero unless there exists a 3-contributing pair $(r,s)$
  such that $((p'rb-psb')^2-(a-3c)^2B^2)/(4Baa')-j$ is divisible by
  $n$.  Equivalently, coefficient $\psi_j$ is zero unless there exists
  an integer $m$ such that $mB(3a'm+a'-3c)-j$ is divisible by $n$, see
  Section \ref{Proof of main5}.
\end{rem}

\begin{rem}\label{rem} 
  The products appearing in the right hand sides of our formulae 
  satisfy some obvious relations.  For
  example for odd $a'$ and $c$, $a'>c$, we have
   \be
  \phi_{a',1,c,1}(q)\prod_{j=1, \ 
  j\neq (c+1)/2}^{(a'-1)/2}\phi_{a',1,2j-1,a'}(q)=1.
    \ee
    Theorems \ref{main} and \ref{main5} can be used to replace
    $\phi_{a',B,c,n}$ and $\psi_{a',B,c,n}$ in this and similar
    formulae via alternating sums of Virasoro characters. That leads
    to identities which involve alternating sums of products of
    Virasoro characters.
\end{rem}

\section{Proofs} \label{proof sec}
\subsection{Proof of Theorem \ref{main}}\label{Proof of main}
The Jacobi triple product identity 
(see for example (2.2.10) in \cite{A1}) reads:
\be
(v,u,u^{-1}v;v)_\infty=
\sum_{j\in\Z}(-1)^ju^jv^{j(j-1)/2}.
\ee
Substituting 
\bean\label{sub}
v=q^{Ba'}, \qquad  u=q^{B(a'+c)/2},
\eean
and changing the
summation index $j$ to $-j$ we obtain the following formula for the
right hand side of \Ref{main for}:
\be
\frac{(q^{B(a'-c)/2},q^{B(a'+c)/2},q^{Ba'};q^{Ba'})_\infty}{(q^n;q^n)_\infty}=
\sum_{j\in\Z}\frac{(-1)^j q^{jB(a'j+c)/2}}{(q^n;q^n)_\infty}\ .
\ee

Substituting further $j=2nk+m$, where $k\in\Z$, $m\in\{0,\dots,2n-1\}$,
we obtain:
\bean\label{inter}
\sum_{j\in\Z}\frac{(-1)^j q^{jB(a'j+c)/2}}{(q^n;q^n)_\infty}=\sum_{m=0}^{2n-1}
(-1)^m q^{mB(a'm+c)/2}
\sum_{k\in \Z} \frac{q^{nkB(2a'nk+2a'm+c)}}{(q^n;q^n)_\infty}\ .
\eean

After substituting Rocha-Caridi formula \Ref{RC} for the Virasoro
characters in the left hand side of formula \Ref{main for}, we obtain
$n$ positive and $n$ negative terms of the form $q^{x_j}\sum_{k\in
  \Z}q^{nk(pp'k+y_j)}/(q^n;q^n)_\infty$ with some $x_j,y_j$.  We claim
that after a linear change of the summation index these terms match
the $2n$ terms in the right hand side of \Ref{inter}.

\begin{lem}\label{adm}
The pair $(r,s)$ is 2-contributing if and only if $0< r< p/b$,
$0< s < p'/b'$, $r$ is odd, $ps \equiv  bc$ $(\on{mod} \ a')$.

There are exactly $n$ 2-contributing pairs. 

If  $(r,s)$ is a 2-contributing pair then
$(p/b-r,p'/b'-s)$ is not a 2-contributing pair.

If $(r,s)$ is a 2-contributing pair then
both $(p'r/b'-ps/b+c)$ and $(p'r/b'+ps/b-c)$ are divisible by $2a'$.
\end{lem}
\begin{proof}
  If $(r,s)$ is a 2-contributing pair then we clearly have that $r$ is
  odd and $ps/b \equiv c$ $(\on{mod} \ a')$. If $r$ is odd and $ps/b
  \equiv c$ $(\on{mod} \ a')$, then clearly $(p'r/b'-pr/b+c)/(2a')$
  and $(p'r/b'+pr/b-c)/(2a')$ are integers.  The sum of these two
  integers equals to $p'r/(a'b')$ which is odd.  Therefore exactly one
  of the numbers $(p'r/b'-pr/b+c)/(2aa')$ and $(p'r/b'+pr/b-c)/(2aa')$
  is an integer and $(r,s)$ is a 2-contributing pair.
  
  Note that $p$ and $a'$ are relatively prime and therefore
  $p(ka'+1),p(ka'+2),\dots,p((k+1)a'-1)$ are all different and
  non-zero modulo $a'$. Therefore exactly one of these numbers has the
  same residue as $bc$ modulo $a'$. It follows that we have
  $p'/(a'b')$ choices for $s$ and similarly we have $p/(ab)$
  independent choices for $r$.  Thus we have $pp'/(aa'B)=n$
  2-contributing pairs.
  
  If $(r,s)$ is a 2-contributing pair then $ps/b\equiv c$ $(\on{mod}\ 
  a')$ and therefore $p(p'/b'-s)/b\equiv -c$ $(\on{mod}\ a')$. Since
  $a'$ is odd, $c$ and $-c$ have different residues and the pair
  $(p/b-r,p'/b'-s)$ is not 2-contributing.
  
  The numbers $p'r/b'\pm ps/b \pm c$ are even integers for all choices
  of pluses and minuses.  Also $ps/b-c$ and $p'r/b'$ are both
  divisible by $a'$. Since $a'$ is odd, the last statement of the
  lemma follows.
\end{proof}

For a 2-contributing pair $(r,s)$, define integers $m_{r,s}$ and 
$\bar m_{r,s}$ as follows. Set $x_{r,s}=1$ if 
$p'r/b'-ps/b+c>0$ and $x_{r,s}=0$ if $p'r/b'-ps/b-c\leq 0$.
Then define
\bean\label{m1}
m_{r,s}&=&2nx_{r,s}-(p'r/b'-ps/b+c)/(2a'), \notag\\ 
\bar m_{r,s}&=&(p'r/b'+ps/b-c)/(2a'). 
\eean
We clearly have $0\le m_{r,s} \le 2n-1$, $0 \le \bar m_{r,s} \le 2n-1$.

\begin{lem} 
The $2n$ numbers $\{m_{r,s},\bar m_{r,s}\}$ are all distinct.
\end{lem}
\begin{proof}
If $m_{r_1,s_1}=m_{r_2,s_2}$ then 
\be
(p'r_1/b'-ps_1/b)-(p'r_2/b'-ps_2/b)=p'(r_1-r_2)/b'-p(s_1-s_2)/b
\ee
is divisible by $4a'n=2pp'/B$.  The divisibility by $p/b$ gives
$r_1=r_2$ and the divisibility by $p'/b'$ gives $s_1=s_2$.

If $\bar m_{r_1,s_1}=\bar m_{r_2,s_2}$ then  
\be
(p'r_1/b'+ps_1/b)-(p'r_2/b'+ps_2/b)=p'(r_1-r_2)/b'+p(s_1-s_2)/b
\ee 
is zero and hence it is divisible by $2pp'/B$. 
Therefore $r_1=r_2$ and $s_1=s_2$. 

If $m_{r_1,s_1}=\bar m_{r_2,s_2}$ then 
\be
(p'r_1/b'-ps_1/b)+(p'r_2/b'+ps_2/b)=p'(r_1+r_2)/b'-p(s_1-s_2)/b
\ee
is divisible by $2pp'/B$. The divisibility by $p/b$ and by $p'/b'$
implies $s_1=s_2$ and $r_1+r_2=p/b$. It leads to a conclusion that
$2pp'/B$ divides $pp'/B$ which is a contradiction.
\end{proof}

\begin{lem}\label{sign} We have 
\be
(-1)^{m_{r,s}}=(-1)^{(p'r/b'-ps/b+c)/{2}}, \qquad {\bar m_{r,s}}=
-(-1)^{(p'r/b'-ps/b+c)/{2}}.
\ee
In particular 
\bea
(-1)^{m_{r,s}}=-(-1)^{\bar m_{r,s}}=1 \qquad {\rm if} \ \ 
(r,s)\in\mc A_1^{(2)},\\
-(-1)^{m_{r,s}}=(-1)^{\bar m_{r,s}}=1 \qquad  {\rm if} \ \ 
(r,s)\in\mc A_2^{(2)}.
\eea
\end{lem}
\begin{proof}
The first equation follows from the definition since $a'$ is odd. 

Since $p/b$ is even and $c$ is odd, we have 
\be
(-1)^{\bar m_{r,s}}= (-1)^{(p'r/b'+ps/b-c)/2}= 
(-1)^{(p'r/b'-ps/b+c)/{2}+(ps/b-c)}
=-(-1)^{(p'r/b'-ps/b+c)/2}.
\ee
The rest of the lemma is obvious.
\end{proof}

Finally, for a 2-contributing pair $(r,s)$ we have
\begin{align*}
(-1)^{\frac{p'r/b'-ps/b+c}{2}}q^{\frac{(p-p')^2-(cB)^2}{8Ba'}}
q^{-n\frac{(p-p')^2}{4pp'}}
\sum_{k\in\Z} q^{\frac{n}{4pp'}(2pp'k+p'rb-psb')^2}=\\
(-1)^{m_{r,s}}\sum_{k\in\Z} q^{\frac{n}{4pp'}
\big(\big(2pp'(-k-x_{r,s})+2pp'x_{r,s}-2a'Bm_{r,s}-cB\big)^2-(cB)^2\big)}=\\
=(-1)^{m_{r,s}} q^{m_{r,s}B(a'm_{r,s}+c)/2}
\sum_{k\in \Z} q^{nkB(2a'nk+2a'm_{r,s}+c)}.
\end{align*}
Similarly:
\begin{align*}
  -(-1)^{\frac{p'r/b'-ps/b+c}{2}}q^{\frac{(p-p')^2-(cB)^2}{8Ba'}}
  q^{-n\frac{(p-p')^2}{4pp'}}\sum_{k\in\Z}
  q^{\frac{n}{4pp'}(2pp'k+p'rb+psb')^2}=\\ =(-1)^{\bar m_{r,s}}
  q^{\bar m_{r,s}B(a'\bar m_{r,s}+c)/2} \sum_{k\in \Z}
  q^{nkB(2a'nk+2a'\bar m_{r,s}+c)}.
\end{align*}
Theorem \ref{main} is proved.

\subsection{Proof of Theorem \ref{main a}}\label{Proof of main a}
The proof of Theorem \ref{main a} is similar to the proof of 
Theorem \ref{main}. The only difference is in signs.

To prove formula \Ref{main for a}, we 
change the substitution \Ref{sub} to
\be
v=-q^{Ba'}, \qquad  u=-q^{B(a'+c)/2},
\ee
and
observe that since $n$ is even, $a'$ is odd,
\be
\frac{m_{r,s}(m_{r,s}-1)}{2}-\frac{\bar m_{r,s}(\bar m_{r,s}-1)}{2}=
\frac{(m_{r,s}-\bar m_{r,s})(m_{r,s}+\bar m_{r,s}-1)}{2}
\ee
has the same parity as
\be
\frac{(p'r/b')(ps/b-c-a')}{2}\ .
\ee
This number is  odd because 
$a'+c$ is even but not divisible by 4, 
$p/b$ is divisible by $4$ and $p'r/b$ is odd.  
This observation replaces Lemma \ref{sign}.

To prove formula \Ref{main for b}, we 
change the substitution \ref{sub} to
\be
v=-q^{Ba'}, \qquad  u=q^{B(a'+c)/2},
\ee
and
observe that since $n$ is even, $a'$ is odd,
\be
\frac{m_{r,s}(m_{r,s}+1)}{2}-\frac{\bar m_{r,s}(\bar m_{r,s}+1)}{2}=
\frac{(m_{r,s}-\bar m_{r,s})(m_{r,s}+\bar m_{r,s}+1)}{2}
\ee
has the same parity as
\be
\frac{(p'r/b')(ps/b-c+a')}{2}\ .
\ee
This number is  odd because 
$a'-c$ is even but not divisible by 4, 
$p/b$ is divisible by $4$ and $p'r/b$ is odd.
This observation replaces Lemma \ref{sign}.

\subsection{Proof of Theorem \ref{pos thm}}\label{Proof of pos}
Let all odd prime divisors of $n$ divide $a'$.

Consider $2n$ numbers $\{Bm(a'm+c)/2,\  m=0,1,\dots,2n-1\}$. We claim that
for each $j\in\{0,\dots,n-1\}$, exactly two of these $2n$ numbers have
residue $j$ modulo $n$.

Consider the following equation for
$x\in\{0,1,\dots ,2n-1\}$:
\bean\label{res eqn}
\frac{Bm_0(a'm_0+c)}{2}\equiv \frac{Bx(a'x+c)}{2}\quad (\on{mod} n).
\eean
To establish our claim, it is sufficient to show that for any 
$m_0\in \{0,1,\dots, 2n-1\}$, 
equation \Ref{res eqn} has exactly two solutions. 

Since $B\perp n$, we cancel $B$ on both sides and obtain that
$(x-m_0)(a'(x+m_0)+c)/2$ is divisible by $n$. Write $n=2^dk$, where
$k$ is odd. Then from our assumptions we have $k\perp (a'(x+m_0)+c)$
and it follows that $k$ divides $x-m_0$. Therefore $x$ has the form
$x=m_0+kl$ for some integer $l$ satisfying $-m_0/k \leq l<(2n-m_0)/k$.

If $x-m_0$ is even then $(a'(x+m_0)+c)$ is odd and it follows that
$x-m_0$ is divisible by $2n$ and therefore $x=m_0$. If $x-m_0$ is odd
then $(a'(x+m_0)+c)$ is divisible by $2^{d+1}$. But since $a'k$ is odd,
the $2^{d+1}$ numbers $\{a'((m_0+kl)+m_0)+c,\ -m_0/k \leq
l<(2n-m_0)/k\}$ all have different residues modulo $2^{d+1}$. Therefore
exactly one of them is divisible by $2^{d+1}$. 

Our claim is proved.

\medskip

If $d>0$ then $B$ is odd and we choose 
$p=2^{d+1}$. In such a case we have $b=1$, $b'=B$, $p'=ka'B$.
If $d=0$ and $B=2^{\tilde d} \tilde k$ with odd $\tilde k$ then we choose 
$p=2^{\tilde d+1}$. In such a case we have $b=2^{\tilde d}$, 
$b'=\tilde k $, $p'=\tilde ka'n$.

Use Theorem \ref{main} to write $\phi(q)$ as a sum of $n$ Virasoro
characters. It follows that for each $j\in\Z_{\geq 0}$ we have exactly
one term in the left hand side of \Ref{main for} which contains $q^j$
and moreover, this term is the same for $j$ and $j+n$. Indeed, each
Virasoro character corresponds to two terms in \Ref{inter} and as we
have shown, exactly two terms in \Ref{inter} contribute to $q^j$ with
a given $j$ modulo $n$. 

The first statement of the theorem is proved.

\medskip

Let now $n$ be an odd prime number and let $B$ be odd. 
(Cf. \cite{A}, proof of Theorem 1.)

Choose $p=2$. We have $p'=a'nB$, $b=1$, $b'=B$.

Consider equation \Ref{res eqn}.
We claim that there are at most $4$ solutions. Indeed
$(x-m_0)(a'x+a'm_0+c)$ is divisible by $n$. Since $n$ is prime then
either $x-m_0$ or $a'x+a'm_0+c$ is divisible by $n$. In each case we
obtain at most two values of $x$.

Therefore for each $j\in\Z_{\geq 0}$ we have at most two terms in the
left  hand side of \Ref{main for}
which contain $q^j$.  But the difference of two
$(2,p')$ Virasoro characters is known to have all coefficients of the
same sign. It follows for example from the fermionic representation of
$(2,p')$ characters used in the Rogers-Ramanujan-Gordon-Andrews
identities, (see (7.3.7) in \cite{A1}).

\subsection{Proof of Theorem \ref{main5}}\label{Proof of main5} 
The proof of Theorem \ref{main5} is similar to that of Theorem
\ref{main}. The main difference is the use of the quintuple product
identity as opposed to the triple product identity.

The quintuple product identity (see \cite{W}) reads:
\be
(v,u,u^{-1}v;v)_\infty (u^2v,u^{-2}v;v^2)_\infty=
\sum_{j\in\Z}(u^{-3j}-u^{3j+1})v^{j(3j+1)/2}.
\ee
Substituting 
\bean\label{sub5}
v=q^{2Ba'}, \qquad  u=q^{Bc}, 
\eean
we obtain the following formula
for the right hand side of \Ref{main5 for}:
\bea
\frac{(q^{Bc},q^{B(2a'-c)},q^{2Ba'};q^{2Ba'})_\infty 
(q^{2B(a'+c)},q^{2B(a'-c)};q^{4Ba'})_\infty}{(q^n;q^n)_\infty}=\\
\sum_{j\in\Z}\frac{q^{jB(3a'j+a'-3c)}-q^{jB(3a'j+a'+3c)+Bc}}{(q^n;q^n)_\infty}.
\eea

Substituting further $j=nk+m$, where $k\in\Z$ and
$m\in\{0,1,\dots,n-1\}$, we obtain:

\bean\label{inter5}
\sum_{j\in\Z}\frac{q^{jB(3a'j+a'-3c)}-q^{jB(3a'j+a'+3c)+Bc}}{(q^n;q^n)_\infty}=
\sum_{m=0}^{n-1}(q^{mB(3a'm+a'-3c)}
\sum_{k\in \Z} \frac{q^{nkB(3a'nk+6a'm+a'-3c)}}{(q^n;q^n)_\infty}- 
\hspace{-30pt}\notag\\
-q^{mB(3a'm+a'+3c)+Bc}
\sum_{k\in \Z} \frac{q^{nkB(3a'nk+6a'm+a'+3c)}}{(q^n;q^n)_\infty})\ .
\eean

After substituting the Rocha-Caridi formula \Ref{RC} for the Virasoro
characters in the left hand side of formula \Ref{main5 for}, we obtain
$n$ positive and $n$ negative terms of the form $q^{x_j}\sum_{k\in
  \Z}q^{nk(pp'k+y_j)}/(q^n;q^n)_\infty$ for some $x_j,y_j$.  We claim
that after a linear change of the summation index these terms match the
$n$ positive and $n$ negative terms in the right hand side of
\Ref{inter5}.

\begin{lem}\label{adm5}
There are exactly $n$ 3-contributing pairs. 

If  $(r,s)$ is a 3-contributing pair then $(p/b-r,p'/b'-s)$ is not a 
3-contributing pair.

The pair $(r,s)$ is a 3-contributing pair of type $1$ if and only if  
$p'r/b'+ps/b-a'-3c$ is divisible by $6a'$.

The pair $(r,s)$ is a 3-contributing pair of type $2$ if and only if
$p'r/b'-ps/b+a'+3c$ is divisible by $6a'$.
\end{lem}
\begin{proof}
  First, consider the case when $p'/(a'b')$ is odd.  

  Then we claim that for any nonnegative integers $k_1,k_2$ such that
  $3k_1<p/b$ and $a'k_2<p'/b'$ there is exactly one 3-contributing
  pair $(r,s)$ such that $3k_1\leq r<3k_1+3$ and $a'k_2\leq s<a'k_2+a'$.
  
  Indeed there is exactly one pair $(r_1,s_1)$ such that
  $3k_1\leq r_1<3k_1+3$, $a'k_2\leq s_1<a'k_2+a'$ and $p'r_1/b'-ps_1/b-a'+3c$
  is divisible by $3a'$.  The numbers $r_1,s_1$ are unique solutions
  (in the specified range) of equations $p'r_1/b'\equiv a'$
  $(\on{mod}\ 3)$ and $ps_1/b\equiv 3c$ $(\on{mod}\ a')$.
  
  Similarly there is exactly one pair $(r_2,s_2)$ such that
  $3k_1\leq r_2<3k_1+3$, $a'k_2\leq s_2<a'k_2+a'$ and $p'r_2/b'+ps_2/b+a'-3c$
  is divisible by $3a'$.
  
  We have $s_1=s_2$. Since $a'\perp 3$ we also have, $r_1\neq r_2$,
  $r_1\neq 3k_1$, $r_2\neq 3k_1$ and therefore $|r_1-r_2|=1$. Recall
  that $p'/(a'b')$ is odd.  It follows that exactly one of the two numbers
  $(p'r_2/b'-ps_2/b-a'+3c)/(3a')$ and $(p'r_1/b'+ps_1/b+a'-3c)/(3a')$ is even 
  and we have exactly one 3-contributing pair.
 
  Now, let $p'/(a'b')$ be even. Then we repeat the same argument.
  However, in this case, the numbers $(p'r_2/b'-ps_2/b-a'+3c)/(3a')$
  and $(p'r_1/b'+ps_1/b+a'-3c)/(3a')$ have the same parity.  But this
  parity is changed when $k_2$ is replaced by $k_2+1$. Therefore for
  half of the possible values of $k_2$ we have two 3-contributing pairs
  and there are no contributing pairs for the other half.

\medskip
 
Let $(r,s)$ be a 3-contributing pair of the first type, that is
$p'r/b'-ps/b-a'+3c$ is divisible by $6a'$. Then
$p'r/b'-ps/b+a'-3c$ and $p'r_1/b'+ps_1/b-a'+3c$ are not divisible
by $6a'$. The number $2pp'/B$ is divisible by $6a'$. It follows that
$(p/b-r,p'/b'-s)$ is not a 3-contributing pair.

The case of a 3-contributing pair of the second type is done similarly.

\medskip

  If $(r,s)$ is a 3-contributing pair of the first type then 
  $2p'r/b'-2a'$ is divisible by $3$. In addition it is clearly
  divisible by $2a'$ and therefore it is divisible by $6a'$.
 
  Similarly, if $(r,s)$ is a 3-contributing pair of the second type 
  then $2p'r/b'+2a'$ is divisible by $6a'$.
 
  The last two statements of the lemma follow.
\end{proof}

If $(r,s)$ is a 3-contributing pair of type $1$, 
we define integers $x_{r,s},\bar x_{r,s}, m_{r,s}, \bar m_{r,s}$ 
by the equality 
\bea
\frac{p'r/b'-ps/b-a'+3c}{6a'}&=&nx_{r,s}+m_{r,s}, \\ 
\frac{p'r/b'+ps/b-a'-3c}{6a'}&=&n\bar x_{r,s}+\bar m_{r,s},
\eea  
and the requirement  $0 \leq  m_{r,s} <n$, $0 \leq  \bar m_{r,s} <n$.

If $(r,s)$ is a 3-contributing pair of type $2$, we define integers
$x_{r,s},\bar x_{r,s}, m_{r,s}, \bar m_{r,s}$ by the equality
\bea
-\frac{p'r/b'+ps/b+a'-3c}{6a'}&=&nx_{r,s}+m_{r,s}, \\
-\frac{p'r/b'-ps/b+a'+3c}{6a'}&=&n\bar x_{r,s}+\bar m_{r,s},
\eea 
and the requirement  $ 0\leq  m_{r,s} <n$, $ 0 \leq \bar m_{r,s} <n$.

\begin{lem} 
The $n$ numbers $\{m_{r,s}\}$ are all distinct. 
The $n$ numbers $\{\bar m_{r,s}\}$ are also all distinct.
\end{lem}
\begin{proof}
If $m_{r_1,s_1}=m_{r_2,s_2}$ then 
\be
(p'r_1/b'-ps_1/b)-(p'r_2/b'-ps_2/b)=p'(r_1-r_2)/b'-p(s_1-s_2)/b
\ee 
is divisible by $2pp'/B$ or
\be
(p'r_1/b'-ps_1/b)+(p'r_2/b'+ps_2/b)=p'(r_1+r_2)/b'-p(s_1-s_2)/b
\ee 
is divisible by $2pp'/B$.

In the former case the divisibility by $p/b$ gives $r_1=r_2$ and the
divisibility by $p'/b'$ gives $s_1=s_2$. In the later case we
similarly obtain $s_1=s_2$ and $r_1+r_2=p/b$.  It leads to a
conclusion that $2pp'/B$ divides $pp'/B$ which is a contradiction.

The case $\bar m_{r_1,s_1}=\bar m_{r_2,s_2}$ is done similarly.
\end{proof}

Finally, for a 3-contributing pair $(r,s)$ of type $1$ we have
\begin{align*}
q^{\frac{(p-p')^2-B^2(a'-3c)^2}{12Ba'}}q^{-n\frac{(p-p')^2}{4pp'}}
\sum_{k\in\Z} q^{\frac{n}{4pp'}(2pp'k+p'rb-psb')^2}=\\
\sum_{k\in\Z} q^{\frac{n}{4pp'}
\big(\big(2pp'(k-\bar x_{r,s})+2pp'x_{r,s}+6a'B m_{r,s}+(a'-3c)B\big)^2-B^2(a'-3c)^2\big)}=\\
=q^{m_{r,s}B(3a'm_{r,s}+a'-3c)}
\sum_{k\in \Z} q^{nkB(3a'nk+6a'm_{r,s}+a'-3c)}.
\end{align*}
Similarly:
\begin{align*}
  q^{\frac{(p-p')^2-B^2(a'-3c)^2}{12Ba'}}q^{-n\frac{(p-p')^2}{4pp'}}
\sum_{k\in\Z} q^{\frac{n}{4pp'}(2pp'k+p'rb+psb')^2}=\\
\sum_{k\in\Z} q^{\frac{n}{4pp'}
\big(\big(2pp'(k-x_{r,s})+2pp'x_{r,s}+6a'B\bar m_{r,s}+(a'+3c)B\big)^2-B^2(a'-3c)^2\big)}=\\
=q^{\bar m_{r,s}B(3a'\bar m_{r,s}+a'+3c)+Bc}
\sum_{k\in \Z} q^{nkB(3a'nk+6a'\bar m_{r,s}+a'+3c)}.
\end{align*}
The computation for a 3-contributing pair of type $2$ is similar.
Theorem \ref{main5} is proved.

\subsection{Proof of Theorem \ref{main5 a}}\label{Proof of main5 a}
The proof of Theorem \ref{main5 a} is similar to the proof of 
Theorem \ref{main5}. The only difference is in signs.

To prove formula \Ref{main5 for a}, we 
change the substitution \Ref{sub5} to
\be
v=q^{2Ba'}, \qquad  u=-q^{Bc},
\ee
and
observe that
$m_{r,s}+\bar m_{r,s}$ has the same parity as $(p'r/b'-a')/a'$. This number 
is clearly odd if $p'/(a'b')$ is even. If $p'/(a'b')$ is odd and $n$ is even then $p/b$ is even, $a'$ is odd and if $c$ is also odd then $r$ is even and 
therefore $(p'r/b'-a')/a'$ is odd.

To prove formula \Ref{main5 for b}, we 
change the substitution \Ref{sub5} to
\be
v=-q^{2Ba'}, \qquad  u=q^{Bc},
\ee
and
observe that
\be
\frac{m_{r,s}(3m_{r,s}+1)}{2}-\frac{\bar m_{r,s}(3\bar m_{r,s}+1)}{2}=
\frac{(m_{r,s}-\bar m_{r,s})(3m_{r,s}+3\bar m_{r,s}+1)}{2}
\ee
has the same parity as
\be
\frac12\ \frac{p'r}{a'b'}\ \frac{(ps/b-3c)}{a'}\ .
\ee
This number is even. Indeed, $n$ is divisible by $4$, hence 
if $p'/(a'b')$ is even then $p'/(a'b')$ is divisible by $4$, and if 
$p'/(a'b')$ is odd then $p/b$ is divisible by $4$.

To prove formula \Ref{main5 for c}, we 
change the substitution \Ref{sub5} to
\be
v=-q^{2Ba'}, \qquad  u=-q^{Bc},
\ee
and
observe that 
\be
\frac{m_{r,s}(3m_{r,s}-5)}{2}-\frac{\bar m_{r,s}(3\bar m_{r,s}-5)}{2}=
\frac{(m_{r,s}-\bar m_{r,s})(3m_{r,s}+3\bar m_{r,s}-5)}{2}
\ee
has the same parity as
\be
\frac12\ \frac{(3p'r-2a'b')}{a'b'}\ \frac{(ps/b-3c)}{a'}\ .
\ee
This number is odd. Indeed if $p/b$ is divisible by $4$ and $c$ is even
then $p'r$ is odd, and if $p'/(a'b')$ is divisible by 
$4$ then $(ps/b-3c)/a'$ is odd.

\subsection{Proof of Theorem \ref{pos thm 5}}\label{Proof of pos 5}
The proof of Theorem \ref{pos thm 5} is similar to that of Theorem
\ref{pos thm}. 

Namely, we write $n=3^dk$ where $k\perp 3$ and show that the $n$ numbers
$Bm(3a'm+a'-3c)$, $m=0,\dots,n-1$, are all different modulo $n$. 
It follows that the $n$ terms of the left hand side of \Ref{main5 for} all
contribute to different coefficients and therefore there is no 
further subtraction.

\end{document}